\documentclass{article}
\usepackage{amsfonts,amsmath}
\mathchardef\mhyphen="2D 
\usepackage[a4paper]{geometry}
\usepackage{enumitem}
\usepackage{url}
\usepackage{microtype}
\usepackage{relsize}

\def\zo/{$0\mkern2mu\mhyphen1$}

\def\nn/{$n \times n$}

\title{The Genius Conjectures (via Bell Polynomials)}
\date{\today} 
\author{Paul Federbush\\
Department of Mathematics\\
University of Michigan\\
Ann Arbor, MI, 48109-1043}

\usepackage{amsthm} 
\newtheorem{thm}{Theorem}[section]
\newtheorem{conj}{Conjecture}
\newtheorem*{conj*}{Conjecture} 
\DeclareMathOperator{\Prob}{Prob}

\newcommand{\stirlingi}{\displaystyle\genfrac{[}{]}{0pt}{}}

\begin{document}

\maketitle
\begin{abstract}
	We present two related conjectures, arising in work on $i$-matchings in
	random $r$-regular bipartite graphs. The conjectures themselves are easily
	stated and involve only basic properties of convergent power series. One
	formulation involves Bell's polynomials. The conjectures name was chosen since
	we earnestly believe only a truly genius mathematician will prove them. We
	advise others not to try. A further belief is that the proof will arise from
	some deep properties of partitions.
\end{abstract}

We write the paper in a form so that the conjectures are reached as early as
possible. Thus Section \ref{sec:1} presents the mathematical setting of the
conjectures, and Section \ref{sec:2} contains the conjectures themselves. One
need read only these two sections to see the task before one. Section
\ref{sec:3} presents a reformulation of the setting using Bell's polynomials
\cite{Wiki}. This section is entirely the work of my colleague David Williams.
One may only read Section \ref{sec:3} and then Section \ref{sec:2} to arrive at
a full treatment of the conjectures using Bell polynomials.

The final section, Section \ref{sec:4}, on background provides the place of the
conjectures in the development of mathematics, how they arose, and what follows
from their truth. Thus, this section contains the contents of usual
introductions and conclusions. I only say now they arose from certain equations
of Mario Pernici \cite{Per} in his treatment of work of Ian Wanless on
$i$-matchings, [3].

\section{Setting}\label{sec:1}

$p$ is a fixed positive integer. Additionally we have:
\begin{itemize}
	\item $x$, $y$ variables
	\item $d_1,d_2,...,d_p$ variables
	\item $u_1=1$
	\item $u_2,...,u_p$ variables
	\item $F_1=1$
	\item $F_2,...,F_p$ variables
\end{itemize}

We use the notation $[x^s]f$ to be the coefficient of $x^s$ in a power series
$f$ in $x$. We require
\begin{equation}\label{eq:1}
	[x^p] e^{\sum_{i=1}^pd_ix^i} = 0
\end{equation}

\begin{thm}
	There are unique functions
	\[
		F_i = F_i(\{u_k\},\{d_k\}) ,	\qquad i=2,...,p
	\]
	such that
	\begin{equation}\label{eq:2}
		[x^p] e^{\sum_{i=1}^p (yu_i+d_i)x^i}
		= [x^p] e^{\sum_{i=1}^p yF_i x^i}
	\end{equation}
	holds
\end{thm}

Notice that \eqref{eq:2} becomes an equality of polynomials in $y$.

The proof is given by extracting a trivial inductive construction of
$F_2,F_3,...,F_p$.

\section{Conjectures}\label{sec:2}

The conjectures concern properties of the functions $F_i$, $i=2,...,p$.

\begin{conj}
	$F_i(\{u_k\},\{d_k\})$ is linear in its dependence on the $\{u_k\}$.
\end{conj}

For the second conjecture we consider varying $p$ over the positive integers, so
$F_i$ becomes a function of $p$: $F_i=F_i(p,\{u_k\},\{d_k\})$.

\begin{conj}\label{conj:2}
	\begin{equation}\label{eq:3}
		F_i = u_i + \sum_{j=1}^{s_i} r_{i,j} m_{i,j} \qquad  i=2,...,p.
\end{equation}
	where $r_{i,j}$ is a rational function of $p$ that goes to zero as $p$ goes
	to infinity, and $m_{i,j}$ is a monomial in the variables
	$\{u_k\}\cup\{d_k\}$.
\end{conj}

For our use of Conjecture \ref{conj:2} a weaker statement would be sufficient.
But the form given has the advantage of containing more of the structure that
might be useful in developing a proof.

\section{Setting via Bell polynomials}\label{sec:3}

In this section we follow David Williams in formulating the setting of Section 1
using Bell polynomials \cite{Wiki}.

$p$ is a fixed positive integer. Additionally we have:

\begin{itemize}
	\item $y$ variable
	\item $d_1,d_2,...,d_p$ variables
	\item $u_1=1$
	\item $u_2,...,u_p$ variables
	\item $F_1=1$
	\item $F_2,...,F_p$ variables
\end{itemize}
We require
\begin{equation}\label{eq:4}
	B_p(1!d_1,...,p!d_p) = 0
\end{equation}
a (scaled) form of \eqref{eq:1}.

\begin{thm}
	There are unique functions
	
	\[
		F_i = F_i(\{u_k\},\{d_k\}) ,	\qquad i=2,...,p
	\]
	such that
	\begin{equation}\label{eq:5}
		B_p [ 1! (yu_1+d_1) , ... , p! (yu_p+d_p) ]
		= B_p [ 1! yF_1 , ... , p! yF_p] .
	\end{equation}
\end{thm}

Equation \eqref{eq:5} is a (scaled) form of \eqref{eq:2}.

The left side of equation \eqref{eq:5} can be written as
\begin{equation}\label{eq:6}
	\sum_{i=0}^{p} \binom{p}{i}
		B_{p-i} [ 1! yu_1 , ... , (p-i)! y u_{p-i} ]
		B_i ( 1!d_1 , ... , i! d_i ) .
\end{equation}

\section{Background}\label{sec:4}

The conjectures of this paper arose in the study of graph positivity which we
now summarize.

We deal with $r$-regular bipartite graphs with $v = 2n$ vertices. We let $m_i$
be the number of $i$-matchings. In 
\cite{BFP}, Butera, Pernici, and I introduced the quantity $d(i)$, in eq.\ $(10)$
therein,
\begin{equation}\label{eq:7}
  d(i) \equiv \ln\biggl( \frac{m_i}{r^i} \biggr)
  - \ln\biggl( \frac{\overline{m}_i}{(v-1)^i} \biggr)
\end{equation}
where $\overline{m}_i$ is the number of $i$-matchings for the complete (not
bipartite complete) graph on the same vertices,
\begin{equation}\label{eq:8}
  \overline{m}_i = \frac{v!}{(v-2i)!\,i!\,2^i}
\end{equation}
We here have changed some of the notation from \cite{BFP} to agree with notation
in \cite{Per}. We then considered $\Delta^k d(i)$ where $\Delta$ is the finite
difference operator, so
\begin{equation}\label{eq:9}
  \Delta d(i) = d(i+1) - d(i)
\end{equation}
A graph was defined to satisfy \textit{graph positivity} if all the meaningful
$\Delta^k d(i)$ were non-negative.
That is
\begin{equation}\label{eq:10}
  \Delta^k d(i) \ge 0
\end{equation}
for $k = 0,\ldots,v$ and $i = 0,\ldots,v-k$.
We made the conjecture, the ``graph positivity conjecture'', supported by some
computer evidence.
\begin{conj*}
	For  fixed $r$, as $n$ goes to infinity the fraction of graphs that satisfy
	graph positivity approaches one.
\end{conj*}
We note some of the impressive results of the numerical study of graph
positivity in \cite{BFP}.
\begin{enumerate}[label=(\arabic*)]
	\item All graphs $v<14$ satisfy graph positivity.
	\item When $r=4$ the first violations occur when $v=22$ in $2$ graphs out of
	the 2806490 graphs with $v=22$.
	\item For $r=3$ the fraction of graphs not satisfying graph positivity
	continuously decreases between $v=14$ and $v=30$. (There is a single
	violation at $v=14$).
\end{enumerate}
We have been working, \cite{Fed1}, to prove the weaker result, weak graph
positivity, the statement
\begin{conj}
	For fixed $r$ and each $i$ and $k$ one has
	\begin{equation}\label{eq:11}
		\Prob \left( \Delta^k d(i) \ge 0 \right) \xrightarrow[n\to\infty]{} 1 .
	\end{equation}
\end{conj}

In our effort toward proving this conjecture a central role is played by
Pernici's work \cite{Per} systematizing results of Wanless, \cite{Wan}. We now
note some definitions from \cite{Per} in slightly modified notation.
\begin{align}
	T_r &= \dfrac{2(r-1)}{2(r-1) - r + r\sqrt{1-4x(r-1)}} .
		\label{eq:12} \\
	u_s(r) &= 2[x^s] T_r .
		\label{eq:13} \\
	M_j &= [x^j] e^{nrx - \sum_{s\ge2} \frac{n u_s(r)}{s} (-x)^s} .
		\label{eq:14} \\
	\left( 1 + \sum_{h=1}^{j-1} \frac{a_h(r,j)}{n^h} \right)
	&= \dfrac{j!}{n^j r^j} M_j .
		\label{eq:15}
\end{align}
These incorporate equations $(3)$, $(10)$ and $(12)$ of \cite{Per}.

We note that $M_j$ may be viewed as the number of $j$-matchings on the
$r$-regular bipartite graph with $2n$ vertices, and having no closed loops, a
non-existent ideal.

In \cite{Per} Pernici presents what I would call a formal derivation of the
following equations:
\begin{align}
	[j^k n^{-h}] \ln \left( 1 + \sum_{s=1}^{j-1} \frac{a_s(r,j)}{n^s} \right)
	&= 0, \qquad k \ge h+2
		\label{eq:16} \\
	[j^{h+1} n^{-h}] \ln \left( 1 + \sum_{s=1}^{j-1} \frac{a_s(r,j)}{n^s} \right)
	&= \frac{1}{(h+1)h} \left( \frac{1}{r^h} - 2 \right)
		\label{eq:17}
\end{align}
and did a significant numerical check of their validity. These equations are
very important to us. As we discuss below we have found a rigorous proof of
\eqref{eq:16}, in fact of a much stronger result than \eqref{eq:16}. As to
\eqref{eq:17}, we think there should be a fussy technical upgrade of Pernici's
formal argument that provides a rigorous proof.

From computer study we came to believe in a stronger form of \eqref{eq:16}.
Namely eq. \eqref{eq:16} holds if one computes the $a_s(r,j)$ instead of using
\eqref{eq:12} and \eqref{eq:13} to compute the quantities $u_s(r)$, using any
values of the $u_s(r)$! Assuming this stronger conjecture I arrived at a
conjecture for Stirling numbers I put on the web, \cite{Fed2}. Robin Chapman has
shown me a proof of this conjecture in a private communication. This conjecture
now proved is presented in the Appendix. From this I was able to prove the
stronger form of \eqref{eq:16}, \cite{Fed3}. This clever response of Robin
Chapman has increased my expectation that someone will prove the Genius
conjectures.

To complete the proof of weak graph positivity, we need a proof of the
conjecture from Section 10 of \cite{Fed1}:
\begin{conj} The Awesome conjecture.
	Let $z_i$ be positive integers. We set
	\begin{equation}\label{eq:18}
		F = \sum_{s\ge0} \dfrac{a_s(r,j)}{n^s}
		+ \sum_i c_i j(j-1)\cdots(j-z_i+1) \dfrac{1}{n^{z_i}r^{z_i}}
			\sum_{s\ge0} \dfrac{a_s(r,j-z_i)}{n^s} .
	\end{equation}
	with $a_0=1$.
	Then we conjecture:
	\begin{align}
		[j^k n^{-h}] \ln (F) &= 0, \qquad k \ge h+2
			\label{eq:19} \\
		[j^{h+1} n^{-h}] \ln (F)
		&= \frac{1}{(h+1)h} \left( \frac{1}{r^h} - 2 \right)
			\label{eq:20}
	\end{align}
\end{conj}
Note this includes \eqref{eq:16} and \eqref{eq:17} above. We will show in
\cite{Fed3} that this Awesome conjecture follows from the validity of
\eqref{eq:17} and the Genius conjectures, and the work of Robin Chapman
\cite{Cha} (leading to the generalization of \eqref{eq:16}). We find the
interrelation of these different conjectures a beautiful context
\bigskip

\noindent 
\begin{Large}
\textbf{APPENDIX Conjecture Proved by Chapman}
\end{Large}
\bigskip

The (unsigned) Stirling numbers of the first kind, $\stirlingi{a}{b}$, are
defined by
\begin{equation}\label{eq:A1}
	x (x+1) \cdots (x+n-1) = \sum_{k=0}^n \stirlingi{n}{k} x^k .
\end{equation}
It is easy to show $\stirlingi{n}{n-w}$ is a polynomial in $n$ of degree $2w$.
So we may naturally define $\stirlingi{x}{x-w}$ for any number $x$ by extending
the domain of the polynomial. We set
\begin{equation}\label{eq:A2}
	P_w(x) \equiv \stirlingi{x}{x-w} .
\end{equation}

Now we give ourself an integer $g\ge2$, an integer $w$, $0\le w\le g-2$, and a
set of $g$ distinct numbers,
\begin{equation}\label{eq:A3}
	S = \{c_1,...,c_g\} .
\end{equation}
We define a \emph{configuration} as a sequence of non-empty subsets of $S$
\begin{equation}\label{eq:A4}
	S_1, S_2, ..., S_r
\end{equation}
that are disjoint with union $S$, i.e.
\begin{equation}\label{eq:A5}
	S_i \ne \emptyset ,	\quad
	S_i \cap  S_j = \emptyset \text{ if } i \ne j ,	\quad
	\bigcup_{i=1}^r S_i = S .
\end{equation}
For a configuration we define
\begin{equation}\label{eq:A6}
	t_i = \sum_{c_k\in S_i}c_k ,	\qquad i=1,...,r .
\end{equation}
A \emph{weighted configuration} is a configuration as above for which each $S_i$
is assigned a non-negative integer, $w_i$, its weight, with the restriction
\begin{equation}\label{eq:A7}
	\sum_{i=1}^r w_i = w .
\end{equation}
Such a weighted configuration has an \emph{evaluation} defined as
\begin{equation}\label{eq:A8}
	(-1)^r \dfrac{1}{r} \prod_i P_{w_i}(t_i) .
\end{equation}

\begin{thm}
	The sum over all distinct weighted configurations of their evaluations is
	zero.
\end{thm}


\begin{thebibliography}{1}
	
  \bibitem{Wiki}
    Bell polynomials, 
    \url{https://en.wikipedia.org/wiki/Bell_polynomials}
   	2020

  \bibitem{Per}
    Pernici, M., \textit{$1/n$ expansion for the number of matchings on regular
    graphs and monomer-dimer entropy,} J. Stat.\ Phys.\ \textbf{168} (2017) 666.

  \bibitem{Wan}
    Wanless, I. M., \textit{Counting Matchings and Tree-like Walks in Regular
    Graphs,} Combinatorics, Probability and Computing \textbf{19} (2010) 463.

  \bibitem{BFP}
    Butera, P., Federbush, P., and Pernici, M., \textit{A positivity property of
    the dimer entropy of graphs,} Physica A \textbf{421} (2015) 208.

  \bibitem{Fed1}
    Federbush, P., \textit{A Near Proof of Weak Graph Positivity, A New Property
    of Regular Random Graphs,} arXiv:1710.00357.

  \bibitem{Fed2}
    Federbush, P., \textit{A Set of Conjectured Identities for Stirling Numbers
    of the First Kind,} arXiv:1808.09264.

  \bibitem{Cha}
   Robin Chapman, private communication.
  
  \bibitem{Fed3}
    Federbush, P., \textit{A Near Proof of Weak Graph Positivity, A New Property
    of Regular Random Graphs, II,} in preparation.

\end{thebibliography}
\end{document}